\documentclass{article}
\usepackage{amsmath, amssymb, latexsym}
\title{The Charney-Davis conjecture for certain subdivisions of spheres}
\author{Andrew Frohmader}

\newtheorem{theorem}{Theorem}[section]
\newtheorem{proposition}[theorem]{Proposition}

\newtheorem{lemma}[theorem]{Lemma}
\newtheorem{definition}[theorem]{Definition}
\newtheorem{conjecture}[theorem]{Conjecture}

\def\proof{\smallskip\noindent {\it Proof: \ }}
\def\endproof{\hfill\ensuremath{\square}\medskip}

\begin{document}

\maketitle

\begin{abstract}
Notions of sesquiconstructible complexes and odd iterated stellar subdivisions are introduced, and some of their basic properties are verified.  The Charney-Davis conjecture is then proven for odd iterated stellar subdivisions of sesquiconstructible balls and spheres.
\end{abstract}

\section{Introduction}

This paper proves the Charney-Davis conjecture for a certain class of simplicial spheres.  We start by discussing our results.  Precise definitions are given in subsequent sections.

Inspired by a longstanding conjecture of Hopf from differential geometry, in 1993, Charney and Davis \cite{charney} conjectured a linear inequality on the face numbers of simplicial spheres that are also flag complexes.  The inequality of their conjecture often fails for spheres that are not flag.

The Charney-Davis conjecture only makes sense in odd dimensions.  In dimension one, it is trivial.  In dimension three, it was proven by Davis and Okun \cite{char3} using heavy topological machinery.  In higher dimensions, it remains open.

For a complex to be flag requires that it not be too ``small".  For example, a simplicial sphere of dimension $d-1$ can have as few as $d+1$ facets, but if it is a flag complex, it must have at least $2^d$ facets.  As such, subdivisions of complexes to make the complex larger (that is, have far more faces of all dimensions) are sensible cases to try.  Taking subdivisions of a complex leaves the complex topologically unchanged, so that a subdivision of a sphere is also a sphere.

Stanley \cite{char1} proved the Charney-Davis conjecture for barycentric subdivisions of shellable spheres, and more generally of S-shellable CW-spheres.  Karu \cite{karu} proved the conjecture for barycentric subdivisions of CW-spheres.  Both of their results were to show that the $cd$-index of the relevant spheres is nonnegative.  That this proves that the barycentric subdivision of such spheres satisfies the Charney-Davis conjecture was observed by Babson.  A barycentric subdivision of a sphere is both a sphere and a flag complex, so these are cases of the conjecture.

Still, a barycentric subdivision of a complex creates far more faces than necessary to make the complex flag, as a barycentric subdivision of a complex of dimension $d-1$ will have $d$! times as many facets as the original complex.  In this paper, we introduce a notion of an odd iterated stellar subdivision, and prove the conjecture for these subdivisions.  The finest subdivision of this type yields the standard barycentric subdivision of a complex.  In general, these are coarser subdivisions than the barycentric subdivision.

We also introduce a notion of sesquiconstructible complexes.  This is a modification of the usual definition of constructible complexes designed to make the proof work.  This is analogous to how Stanley introduced the notion of S-shellability to make his proof work in \cite{char1}.  It is possible that the notions of constructible and sesquiconstructible could coincide.

It is not immediately obvious from the definition that sesquiconstructible complexes other than simplexes actually exist, however.  We show that all shellable complexes are sesquiconstructible.  Furthermore, there are complexes that are sesquiconstructible but not shellable, so this is a strictly weaker requirement.

Our proof applies to odd iterated stellar subdivisions of sesquiconstructible balls and spheres.  These constructions are not necessarily flag complexes.  While the Charney-Davis conjecture only applies to spheres, our proof necessarily proves the conjecture for subdivisions of balls, as well.  This is not a big leap, as our construction necessarily has every interior edge contain at least one interior vertex, and if the Charney-Davis conjecture failed for such a ball, a sphere constructed by attaching two copies of the ball along their boundaries would yield a counterexample.

In Section~2, we discuss shellable and constructible complexes.  We introduce sesquiconstructible complexes, and prove some of their basic properties.  In Section~3, we review stellar subdivisions.  We then introduce the concept of odd iterated stellar subdivisions.  In Section~4, we discuss the Charney-Davis conjecture.  Finally, in Section~5, we give our proof.  The proof is combinatorial in nature, and does not require a background in topology or algebraic geometry.

\section{Shellable, constructible, and sesquiconstructible}

Recall that a \textit{simplicial complex} $\Delta$ on a vertex set $W$ is a collection of subsets of $W$ such that (i) for every $v \in W$, $\{v\} \in \Delta$ and (ii) for every $B \in \Delta$, if $A \subset B$, then $A \in \Delta$.  The elements of $\Delta$ are called \textit{faces}.  A face on $i$ vertices is said to have \textit{dimension} $i-1$, while the dimension of a complex is maximum dimension of a face of the complex.  A simplicial complex in which all maximal faces (under inclusion) are of the same dimension is called \textit{pure}.  The maximal faces are called \textit{facets}.  Faces of dimension one less than facets are called \textit{ridges}.

We can define the \textit{boundary} of a pure complex $\partial \Delta$ by declaring a ridge to be in the boundary if it is contained in only one facet.  A face is contained in $\partial \Delta$ if it is contained in a boundary ridge of $\Delta$.

A simplicial complex is called a (simplicial) \textit{ball} or \textit{sphere} if its geometric realization is homeomorphic to a ball or sphere, respectively.  Every ridge of a ball or sphere is contained in at most two facets.  The boundary of a ball is a sphere of dimension one less than the ball, and the boundary of a sphere is $\emptyset$.

We are primarily concerned with complexes that can be built up from smaller complexes in a nice manner.  Two important notions of such nice complexes are shellable and constructible complexes.

\begin{definition}
\textup{A pure simplicial complex $\Delta$ of dimension $d-1$ is \textit{shellable} if either $\Delta$ is a simplex or $\Delta = U \cup V$, where}
\begin{enumerate}
\item \textup{$U$ is a shellable complex of dimension $d-1$,}
\item \textup{$V$ is a simplex of dimension $d-1$, and}
\item \textup{$U \cap V$ is a shellable complex of dimension $d-2$.}
\end{enumerate}
\end{definition}

Intuitively, a shellable complex can be built up by adding one facet at a time and be a nice complex each step of the way.  The order in which facets are added is called a shelling order.  Shellable complexes have been extensively studied, with results such as that the boundary of any polytope is shellable \cite{polyshell}.

\begin{definition}
\textup{A pure simplicial complex $\Delta$ of dimension $d-1$ is \textit{constructible} if either $\Delta$ is a simplex or $\Delta = U \cup V$, where}
\begin{enumerate}
\item \textup{$U$ and $V$ are constructible complexes of dimension $d-1$, and}
\item \textup{$U \cap V$ is a constructible complex of dimension $d-2$.}
\end{enumerate}
\end{definition}

Constructible complexes may have first been explicitly defined by Hochster \cite{hochster}, though the concept had been around before, as in Zeeman's book \cite{zeeman}.

Note that requiring $V$ to be a simplex yields the definition of shellable complexes, so all shellable complexes are constructible.  Constructible complexes are a more general notion that allows some facets to be attached to each other before needing to being added to the rest of the complex, and there are complexes that are constructible but not shellable.  Hachimori gives some examples in \cite{hachimori}.

For our purposes in this paper, we introduce a new concept that is closely related to that of constructible complexes.

\begin{definition}
\textup{A pure simplicial complex $\Delta$ of dimension $d-1$ is \textit{sesquiconstructible} if either $\Delta$ is a simplex or $\Delta = U \cup V$, where}
\begin{enumerate}
\item \textup{$U$ and $V$ are sesquiconstructible complexes of dimension $d-1$,}
\item \textup{$U \cap V$ is a sesquiconstructible complex of dimension $d-2$, and}
\item \textup{either $\partial(U \cap V) = \emptyset$ or $\partial(U \cap V)$ is a sesquiconstructible complex of dimension $d-3$.}
\end{enumerate}
\end{definition}

It is clear from the definition that every sesquiconstructible complex is constructible.  It is unclear whether the converse is true, as being sesquiconstructible imposes an extra condition on how complexes can be assembled.

More problematic is that it is not immediately obvious that sesquiconstructible complexes actually exist, apart from simplexes and low dimensional examples.  We address this issue by showing that all shellable complexes are sesquiconstructible.  First we need to cite a theorem and prove a lemma.  For a proof of this next theorem, see \cite[Theorem 11.4]{bjorner}.

\begin{theorem} \label{ridgeball}
Let $\Delta$ be a constructible simplicial complex in which every ridge is contained in at most two facets.  Then $\Delta$ is either a ball or a sphere.
\end{theorem}
Since this theorem applies to all constructible complexes, it likewise applies to sesquiconstructible complexes and shellable complexes.

\begin{lemma} \label{sesquishell}
Let $\Delta$ be a simplex of dimension $d-1$ and $C \subset \partial \Delta$ be a pure complex of dimension $d-2$.  Then either $\partial C = \emptyset$ or else $\partial C$ is a shellable complex of dimension $d-3$.
\end{lemma}

\proof  Since $\partial \Delta$ is a sphere, every ridge is contained in exactly two facets.  In particular, every ridge of $C$ is contained in exactly two facets of $\partial \Delta$, and hence at most two facets of $C$.  Any ordering of facets of $C$ is a shelling order, so $C$ is shellable, and by Theorem~\ref{ridgeball}, $C$ is either a ball or a sphere.

If $C$ is a sphere, then $\partial C = \emptyset$.  If $C$ is a ball, then $\partial C$ is a sphere of dimension $d-3$.  Since $\Delta$ has $d$ vertices, $\partial C$ has at most $d$ vertices.  Thus, $\partial C$ is a sphere of dimension $d-3$ on at most $d$ vertices.  Every such sphere is isomorphic to the boundary of some simplicial polytope (see \cite[p. 88, ex. \#4]{ewald}), and hence shellable.  \endproof

It is also possible to give a more self-contained proof of the previous lemma by constructing a shelling order.

\begin{proposition}  \label{constshell}
Every shellable complex is sesquiconstructible.
\end{proposition}

\proof  We use induction on the dimension.  For the base case, the only simplicial complex of dimension -1 is $\{ \emptyset \}$.  This is a simplex, and hence sesquiconstructible.  For the inductive step, we proceed by induction on the number of facets.  For the base case in this second induction, a complex with only one facet is a simplex, and hence sesquiconstructible.

For the inductive step, suppose that $\Delta$ is a shellable complex of dimension $d-1$ with $n$ facets.  Since $\Delta$ is shellable, we can write $\Delta = U \cup V$ with $U$ shellable, $V$ a facet, and $U \cap V$ shellable of dimension $d-2$.  By the second inductive hypothesis, $U$ and $V$ are sesquiconstructible.  By the first inductive hypothesis, so is $U \cap V$.  By Lemma~\ref{sesquishell}, either $\partial (U \cap V) = \emptyset$ or $\partial (U \cap V)$ is shellable of dimension $d-3$, and hence sesquiconstructible by the first inductive hypothesis.  Either way, $\Delta$ is sesquiconstructible.  \endproof

We can also show that for balls and spheres of dimension at most four, the notions of constructible and sesquiconstructible coincide.

\begin{lemma} \label{ucapv}
Let $\Delta$ be a constructible ball or sphere with at least two facets, and let $\Delta = U \cup V$ as in the definition of constructible.  Then $U \cap V$ is either a ball or a sphere.
\end{lemma}

\proof  Since $\Delta$ is a ball or a sphere, every ridge of $\Delta$ is contained in at most two facets.  Every ridge of $U$ is then contained in at most two facets of $\Delta$, and hence at most two facets of $U$.  This and the assumption that $U$ is constructible give that $U$ is either a ball or a sphere by Theorem~\ref{ridgeball}.

We claim that $U \cap V \subset \partial U$.  If not, then there is a face $G \in (U \cap V)$ with $G \not \in \partial U$.  Since $U \cap V$ is pure, there is a facet $F \in (U \cap V)$ with $G \subset F$; $F$ is a ridge of $U$ and a ridge of $V$.  That $G \not \in \partial U$ gives $F \not \in \partial U$, and so $F$ is contained in more than one facet of $U$.  Since $V$ is pure, $F$ is also in a facet of $V$.  $U$ and $V$ have no common facets, for otherwise $U \cap V$ would be of the same dimension as $\Delta$.  Thus, $F$ is a ridge contained in at least three facets of $\Delta$, so $\Delta$ is not a ball or sphere, a contradiction.

If $U$ is a sphere, $\partial U = \emptyset$, which contradicts the assumption that $U \cap V$ is of dimension one less than $\Delta$ unless dim $\Delta \leq 0$, in which case the lemma trivially holds.  Otherwise $U$ is a ball, and so $\partial U$ is a sphere.  Every ridge of $\partial U$ is thus contained in two facets.  Since $(U \cap V) \subset \partial U$ and dim $(U \cap V) = $ dim $\partial U$, every ridge of $U \cap V$ is contained in at most two facets of $\partial U$, and thus at most two facets of $U \cap V$.  Lemma~\ref{ridgeball} then gives that  $U \cap V$ is either a ball or a sphere because $U \cap V$ is constructible.  \endproof

\begin{proposition} \label{consesqui4}
If $\Delta$ is a constructible ball or sphere of dimension at most four, then $\Delta$ is sesquiconstructible.
\end{proposition}

\proof  We use induction on the dimension of $\Delta$.  For the base case, if $\Delta$ is a simplex, it is sesquiconstructible.  Otherwise, if $U$ and $V$ are as in the definition of constructible and dim $\Delta \leq 0$, then dim $\partial (U \cap V) \leq -2$, so $\partial (U \cap V) = \emptyset$, and the additional condition of sesquiconstructible complexes is trivial.  For the inductive step, we use induction on the number of facets of $\Delta$.  For the base case, if $\Delta$ is a simplex, it is sesquiconstructible by definition.

For the inductive step, since $\Delta$ is constructible, $\Delta = U \cup V$, with $U$ and $V$ constructible, and $U \cap V$ constructible of dimension one less than $\Delta$.  By the second inductive hypothesis, $U$ and $V$ are sesquiconstructible, as they have fewer facets than $\Delta$.  By the first inductive hypothesis, $U \cap V$ is sesquiconstructible, as it is of dimension less than $\Delta$.

By Lemma~\ref{ucapv}, $U \cap V$ is either a ball or a sphere. If $U \cap V$ is a sphere, then $\partial (U \cap V) = \emptyset$, so $\Delta$ is sesquiconstructible.  If $U \cap V$ is a ball, then $\partial (U \cap V)$ is a sphere of dimension at most two.  All spheres of dimension at most two are polytopal (see \cite[p. 85-86]{ewald}), and hence shellable, and therefore sesquiconstructible by Lemma~\ref{sesquishell}.  In this case also, $\Delta$ is sesquiconstructible.  \endproof

As shown by Hachimori \cite{hachimori}, there are balls of dimension 3 that are constructible but not shellable.  By Proposition~\ref{consesqui4}, such balls are sesquiconstructible but not shellable.  Hence, sesquiconstructible is a weaker requirement than shellable.

\section{Stellar subdivisions}

The other critical concept that we need is that of a stellar subdivision, which breaks up faces into smaller faces while leaving the complex topologically unchanged.  This is a standard concept; see, for example, \cite[part II, p. 70, def. 2.1]{ewald}.

\begin{definition}
\textup{A \textit{stellar subdivision} of a simplicial complex $\Delta$ at a non-empty face $F \in \Delta$ is obtained by picking a point $p$ in the (relative) interior of $F$.  For each (closed) face $G \in \Delta$ with $p \in G$, we replace $G$ by a cone at $p$ over all faces of $G$ that do not contain $p$.}
\end{definition}

For example, if we start with the triangle on the left in the next picture, we can get the other pictures by subdividing at an edge or at the facet.  The chosen point $p$ is labeled in both pictures.

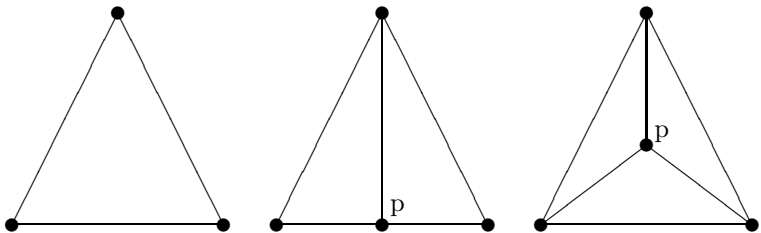
\begin{picture}(300, 100)
\put (10, 10){\circle*{5}}
\put (50, 90){\circle*{5}}
\put (90, 10){\circle*{5}}
\put (10, 10){\line(1, 0){80}}
\put (10, 10){\line(1, 2){40}}
\put (90, 10){\line(-1, 2){40}}

\put (110, 10){\circle*{5}}
\put (150, 90){\circle*{5}}
\put (190, 10){\circle*{5}}
\put (150, 10){\circle*{5}}
\put (110, 10){\line(1, 0){80}}
\put (110, 10){\line(1, 2){40}}
\put (190, 10){\line(-1, 2){40}}
\put (150, 10){\line(0, 1){80}}
\put (153, 15){p}

\put (210, 10){\circle*{5}}
\put (250, 90){\circle*{5}}
\put (290, 10){\circle*{5}}
\put (250, 40){\circle*{5}}
\put (210, 10){\line(1, 0){80}}
\put (210, 10){\line(1, 2){40}}
\put (290, 10){\line(-1, 2){40}}
\put (250, 40){\line(0, 1){50}}
\put (250, 40){\line(4, -3){40}}
\put (250, 40){\line(-4, -3){40}}
\put (253, 43){p}

\end{picture}

We will need couple of lemmas on what happens when taking a stellar subdivision of a face of a complex.

\begin{lemma}  \label{stelsimp}
Let $\Delta$ be a simplex and let $F \in \Delta$ be a non-empty face. Then a stellar subdivision of $\Delta$ at $F$ is sesquiconstructible.
\end{lemma}

\proof  $\partial F$ is the boundary of a simplex, and thus shellable, as every order of its facets is a shelling order.  A cone over a shellable complex is shellable, as a shelling order of the facets before coning is still a shelling order of the cones over what were previously facets.  The stellar subdivision of $F$ at $\Delta$ is obtained by coning over $\partial F$ $(|\Delta| - |F| + 1)$ times, and is thus shellable.  The lemma then follows by Proposition~\ref{constshell}.  \endproof

\begin{lemma}
Let $\Delta$ be a sesquiconstructible complex of dimension $d-1$ and let $F \in \Delta$ be a non-empty face. Then a stellar subdivision of $\Delta$ at $F$ is sesquiconstructible.
\end{lemma}

\proof  We use induction on $d$.  In the base case of $d = 1$, the only non-empty faces are vertices.  A stellar subdivision at a vertex leaves the complex unchanged.  For the inductive step, we use induction on the number of facets of $\Delta$.  In the base case of only one facet, $\Delta$ is a simplex, and the lemma holds by Lemma~\ref{stelsimp}.

To prove the inductive step in this second induction, if $\Delta$ has at least two facets and is sesquiconstructible, then $\Delta = U \cup V$, where $U$ and $V$ are both sesquiconstructible of dimension $d-1$, $U \cap V$ is sesquiconstructible of dimension $d-2$, and $\partial (U \cap V)$ is sesquiconstructible of dimension $d-3$.

By abuse of language, say that if $F \not \in U$, the stellar subdivision of $U$ at $F$ is $U$.  By the inductive hypothesis on the number of facets, the stellar subdivisions of $U$ at $F$ and of $V$ at $F$ are sesquiconstructible of dimension $d-1$.  By the inductive hypothesis on the dimension, the stellar subdivisions of $U \cap V$ and $\partial (U \cap V)$ at $F$ are sesquiconstructible of dimensions $d-2$ and $d-3$, respectively.  Thus, all the conditions for the stellar subdivision of $\Delta$ at $F$ to be sesquiconstructible are satisfied.  \endproof

With this last lemma, if we begin with a sesquiconstructible complex, we can take repeated stellar subdivisions of faces of it in whatever order we like and be guaranteed that the end result will remain sesquiconstructible.

\begin{definition}
\textup{Let $\Delta$ be a simplicial complex.  Let $F_1, F_2, \dots, F_s$ be a subset of non-empty faces of $\Delta$ such that}
\begin{enumerate}
\item \textup{all non-empty faces of odd dimension are in the list and}
\item \textup{if dim $F_i > $ dim $F_j$, then $i < j$.}
\end{enumerate}
\textup{Define complexes $\Delta_0, \Delta_1, \dots, \Delta_s$ such that $\Delta = \Delta_0$ and for all $i > 0$, $\Delta_i$ is the stellar subdivision of $\Delta_{i-1}$ at $F_i$.  We say that $\Delta_s$ is an \textit{odd iterated stellar subdivision} of $\Delta$.}
\end{definition}

Taking a stellar subdivision of a complex at a face can only break up that face and other faces of strictly higher dimension.  As such, a face $F_i$ will be intact in $\Delta_{i-1}$ when it comes time to take a stellar subdivision at $F_i$, so this is well-defined.

The complex $\Delta_s$ depends both on the faces of $\Delta$ chosen and on the relative order of faces of the same dimension in the list.  Choosing the same faces while changing their order can sometimes yield complexes that are not isomorphic to each other.  For example, if we start with a complex consisting of two triangles with a common edge and then do an odd iterated stellar subdivision that divides all edges, we can get either of the complexes shown below, among others.

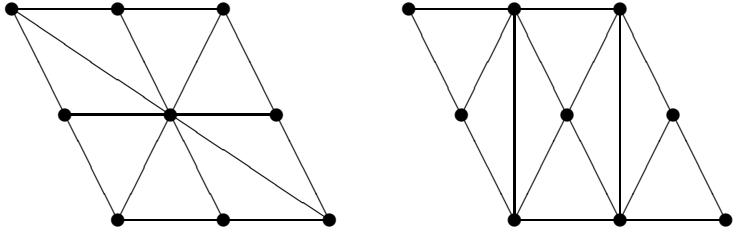
\begin{picture}(300, 100)
\put (10, 90){\circle*{5}}
\put (50, 90){\circle*{5}}
\put (90, 90){\circle*{5}}
\put (30, 50){\circle*{5}}
\put (70, 50){\circle*{5}}
\put (110, 50){\circle*{5}}
\put (50, 10){\circle*{5}}
\put (90, 10){\circle*{5}}
\put (130, 10){\circle*{5}}
\put (10, 90){\line(1, 0){80}}
\put (10, 90){\line(1, -2){40}}
\put (90, 90){\line(1, -2){40}}
\put (50, 10){\line(1, 0){80}}
\put (50, 10){\line(1, 2){40}}
\put (10, 90){\line(3, -2){120}}
\put (50, 90){\line(1, -2){40}}
\put (30, 50){\line(1, 0){80}}

\put (160, 90){\circle*{5}}
\put (200, 90){\circle*{5}}
\put (240, 90){\circle*{5}}
\put (180, 50){\circle*{5}}
\put (220, 50){\circle*{5}}
\put (260, 50){\circle*{5}}
\put (200, 10){\circle*{5}}
\put (240, 10){\circle*{5}}
\put (280, 10){\circle*{5}}
\put (160, 90){\line(1, 0){80}}
\put (160, 90){\line(1, -2){40}}
\put (240, 90){\line(1, -2){40}}
\put (200, 10){\line(1, 0){80}}
\put (200, 10){\line(1, 2){40}}
\put (200, 90){\line(0, -1){80}}
\put (240, 90){\line(0, -1){80}}
\put (180, 50){\line(1, 2){20}}
\put (240, 10){\line(1, 2){20}}
\put (200, 90){\line(1, -2){40}}

\end{picture}

If $F_1, \dots, F_s$ contains all of the faces of $\Delta$ of dimension at least one, then $\Delta_s$ is the barycentric subdivision of $\Delta$.  In general, $\Delta_s$ is a coarser subdivision of $\Delta$.

\section{The Charney-Davis conjecture}

A simplicial complex is a \textit{flag complex} if all of its
minimal non-faces are two element sets. Equivalently, if all of the edges of a potential face of a flag complex are in the
complex, then that face must also be in the complex.

Flag complexes are closely related to graphs.  Given a graph $G$, define its \textit{clique complex} $C = C(G)$ as the simplicial complex whose vertex set is the vertex set of $G$, and whose faces are the cliques of $G$.  The clique complex of any graph is itself a flag complex, as for a subset of vertices of a graph not to form a clique, two of them must not form an edge. Conversely, any flag complex is the clique complex of its 1-skeleton.

Given a simplicial complex, we can count its faces of various dimensions.

\begin{definition}
\textup{The \textit{$i$-th face number} of a simplicial complex $\Delta$, denoted $f_i(\Delta)$ is the number of faces in $\Delta$ of dimension $i$. These are also called \textit{$i$-faces} of $\Delta$. If dim $\Delta = d-1$, the \textit{face vector} of $\Delta$ is the vector
$$f(\Delta) = (f_{-1}(\Delta), f_0(\Delta), \dots , f_{d-1}(\Delta)).$$}
\end{definition}

It is sometimes more convenient to work with a certain linear transform of the face vector.

\begin{definition}
\textup{Define the \textit{h-vector} of a complex $\Delta$ of dimension $d-1$ by $h(\Delta) = (h_0(\Delta), \dots, h_d(\Delta))$, where
$$h_i(\Delta) = \sum_{k=0}^i(-1)^{i-k}{d-k \choose d-i}f_{k-1}(\Delta), \quad 0 \leq i \leq d.$$}
\end{definition}

The h-vector and the face vector of a complex contain exactly the same information, as it is possible to uniquely determine either one from the other.

\begin{theorem}[Dehn-Sommerville relations \cite{dehns}]
If $\Delta$ is a sphere (or more generally, an Eulerian complex) of dimension $d-1$, then $h_i(\Delta) = h_{d-i}(\Delta)$ for $0 \leq i < {d \over 2}$.
\end{theorem}

This theorem says that for a sphere, knowing the first half of the h-vector is enough to determine the rest of it.  Likewise, knowing the first half of the face vector of a sphere uniquely determines the rest of it.

\begin{conjecture}[Charney-Davis \cite{charney}]
Let $\Delta$ be a flag complex of odd dimension $2d-1$.  Suppose further that $\Delta$ is a sphere (or more generally, a Gorenstein* complex).  Then
$$(-1)^d\sum_{i=0}^{2d}(-1)^i h_i(\Delta)\geq 0.$$
\end{conjecture}

Note that such a conjecture is trivial in even dimension, as an alternating sum of the h-numbers of an even dimensional sphere is zero by the Dehn-Sommerville relations.

For our purposes, it is more convenient to express the Charney-Davis conjecture in terms of face numbers rather than h-numbers.  It is straightforward to show that for any simplicial complex $\Delta$ of dimension $n-1$, $$\sum_{i=0}^{n}(-1)^i h_i(\Delta) = 2^n\sum_{i=0}^{n}\bigg({-1 \over 2}\bigg)^i f_{i-1}(\Delta).$$
As such, the form we use is to define the \textit{Charney-Davis quantity}
$$\kappa(\Delta) = \sum_{i=0}^{2d}\bigg({-1 \over 2}\bigg)^i f_{i-1}(\Delta).$$
The Charney-Davis conjecture then states that $\kappa(\Delta) \geq 0$ if $d$ is even and $\kappa(\Delta) \leq 0$ if $d$ is odd.  If $\Delta$ has even dimension $2d$, we can likewise define
$$\kappa(\Delta) = \sum_{i=0}^{2d+1}\bigg({-1 \over 2}\bigg)^i f_{i-1}(\Delta).$$

\section{The main theorem}

Now we can state the main theorem that we wish to prove.

\begin{theorem}  \label{sesquisphere}
Let $\Delta$ be a sesquiconstructible ball or sphere of odd dimension $2d-1$, and let $S$ be an odd iterated stellar subdivision of $\Delta$.  Then $\kappa(S) \geq 0$ if $d$ is even and $\kappa(S) \leq 0$ if $d$ is odd.  Therefore, the Charney-Davis conjecture holds for this class of complexes.
\end{theorem}

Before we can prove this theorem, we need a lemma.

\begin{lemma}  \label{ballhalf}
Let $\Delta$ be a ball or sphere of even dimension $2d$.  Then $2 \kappa(\Delta) = \kappa(\partial \Delta).$
\end{lemma}

\proof  If $\Delta$ is a ball, create a sphere $S$ by taking the union of $\Delta$ and a cone over $\partial \Delta$.  If $\Delta$ is a sphere, take $S = \Delta$.  Since $S$ is a sphere of even dimension, $\kappa(S) = 0$.

We wish to compute $\kappa(S) - \kappa(\Delta)$.  Each face $F \in S$ with $F \not \in \Delta$ was created by coning over $\partial \Delta$.  It thus corresponds to a face of $\partial \Delta$ of dimension (dim $F) - 1$.  Conversely, each face of $\partial \Delta$ corresponds to a face of $S - \Delta$ of dimension one greater.

The face of $\partial \Delta$ contributes $\big({-1 \over 2}\big)^{\textup{dim}F}$ to $\kappa(\partial \Delta)$, and $F$ contributes $\big({-1 \over 2}\big)^{\textup{dim}F + 1}$ to $\kappa(S) - \kappa(\Delta)$.  Hence, each contribution to $\kappa(\partial \Delta)$ is paired with a contribution of ${-1 \over 2}$ times it to $\kappa(S) - \kappa(\Delta)$.  Therefore, $\kappa(S) - \kappa(\Delta) = {-1 \over 2}\kappa(\partial \Delta)$.  Plugging in $\kappa(S) = 0$ and multiplying by $-2$ yields completes the proof.  \endproof

\smallskip\noindent {\it Proof of Theorem~\ref{sesquisphere}: \ }  We use induction on $d$.  For the base case, if $d = 0$, then $\Delta = \{ \emptyset \}$, so $\kappa(\Delta) = 1 \geq 0$.

For the inductive step, we use induction on the number of facets of $\Delta$.  For the base case, if $\Delta$ is a simplex, then $\Delta_1$ is a cone over the $2d-2$ dimensional sphere $\partial \Delta$.  Subsequent stellar subdivisions subdivide this sphere, but $\Delta_s$ is still a cone over a $2d-2$ dimensional sphere.  As such, by Dehn-Sommerville, $\kappa(\partial \Delta_s) = 0$.  Lemma~\ref{ballhalf} then gives that $\kappa(\Delta_s) = 0$.

For the inductive step in the second induction, if $\Delta$ is sesquiconstructible with more than one facet, then $\Delta = U \cup V$ as in the definition of sesquiconstructible.  By Lemma~\ref{ucapv}, $U \cap V$ is either a ball or a sphere and so  Lemma~\ref{ballhalf} gives that $\kappa(U \cap V) = {1 \over 2}\kappa(\partial (U \cap V))$.

By inclusion-exclusion, for any $i$, $f_i(\Delta) = f_i(U) + f_i(V) - f_i(U \cap V)$.  Since $\kappa(\Delta)$ is a linear combination of $f_i(\Delta)$s, $\kappa(\Delta) = \kappa(U) + \kappa(V) - \kappa(U \cap V) = \kappa(U) + \kappa(V) - {1 \over 2}\kappa(\partial (U \cap V))$.

Note that an odd iterated stellar subdivision of a complex induces an odd iterated stellar subdivision of any subcomplex of it.  As such, the inductive hypothesis (for induction on $d$) applies to $\partial (U \cap V)$.  If $d$ is even, then by the inductive hypotheses, $\kappa(U) \geq 0$, $\kappa(V) \geq 0$, and $\kappa(\partial (U \cap V)) \leq 0$, so $\kappa(\Delta) \geq 0$.  If $d$ is odd, then $\kappa(U) \leq 0$, $\kappa(V) \leq 0$, and $\kappa(\partial (U \cap V)) \geq 0$, so $\kappa(\Delta) \leq 0$.  Either way, this completes the inductive step.  \endproof

Gal's conjecture \cite{gal} generalizes the Charney-Davis conjecture for flag spheres.  He conjectures that all of the coefficients of a certain polynomial associated with a sphere are positive if the sphere is flag, and one coefficient is the usual Charney-Davis quantity.  It is unclear whether his conjecture holds for odd iterated stellar subdivisions of sesquiconstructible spheres.

\textit{Acknowledgements.}  I would like to thank Eran Nevo for suggesting a simpler proof of Lemma~\ref{ballhalf}, as well as for pointing out that the odd iterated stellar subdivisions are not always flag complexes.

\end{document}